\documentclass[11pt,twoside]{amsart}

\usepackage{enumerate}
\usepackage[all]{xy}
\usepackage[hyperindex=true,plainpages=false,colorlinks=false,pdfpagelabels]{hyperref}
\usepackage{verbatim}
\usepackage{graphicx,color}
\usepackage{subfigure}
\theoremstyle{plain}
\newtheorem{The}{Theorem}%[section]
\newtheorem*{The*}{Theorem}
\newtheorem{Pro}{Proposition}[section]
\newtheorem{Lem}{Lemma}[section]

\newtheorem*{Cor*}{Corollary}

\theoremstyle{definition}
\newtheorem{Rem}{Remark}[section]

\newtheorem*{Rem*}{Remark}
\numberwithin{equation}{section}

\DeclareMathOperator{\SL}{SL}

\DeclareMathOperator{\Tr}{tr}             % trace

\DeclareMathOperator{\vol}{vol}

\renewcommand{\Re}{\operatorname{Re}}

\newcommand{\dvector}[1]{{\left(\begin{matrix}#1\end{matrix}\right)}}
\DeclareMathOperator{\dbar}{\bar\partial}

\newcommand{\R}{\mathbb{R}}
\newcommand{\C}{\mathbb{C}}
\newcommand{\N}{\mathbb{N}}
\newcommand{\Z}{\mathbb{Z}}

\begin{document}

\title[The asymptotic behavior of the  monodromy representations]{The asymptotic behavior of the  monodromy representations of the associated families  of compact CMC surfaces}

\author{Sebastian Heller}

\address{Sebastian Heller \\
  Institut f\"ur Mathematik\\
 Universit{\"a}t T\"ubingen\\ Auf der Morgenstelle
10\\ 72076 T\"ubingen\\ Germany
 }

\email{heller@mathematik.uni-tuebingen.de}

%\subjclass{53A10,53C42,53C43,14H60}

\date{\today}

\thanks{The author would like to thank Aaron Gerding for many interesting discussions. Author supported by DFG HE 6829/1-2}

\begin{abstract} 
Constant mean curvature (CMC) surfaces in space forms can be described by their associated  $\mathbb C^*$-family of flat $\SL(2,\C)$-connections
 $\nabla^\lambda$.
In this paper we consider the asymptotic behavior (for $\lambda\to0$) of the gauge equivalence classes of $\nabla^\lambda$ for compact CMC surfaces of genus $g\geq2.$ 
We
 prove (under the assumption of simple umbilics) that  the asymptotic behavior  of the traces of the monodromy representation of $\nabla^{\lambda}$ determines the conformal type as well as the Hopf differential
 locally in the Teichm\"uller space.
  \end{abstract}

\maketitle

\setcounter{tocdepth}{1}
%\tableofcontents

%%%%%%%%%%%%%%%%%%%%%%%%%%%%%%%%%%%%%%%%%%%%%%%%%%%%%%%%%%%%%%%%%%%%%%%%%%%%%%
%           Introduction                                                     %
%%%%%%%%%%%%%%%%%%%%%%%%%%%%%%%%%%%%%%%%%%%%%%%%%%%%%%%%%%%%%%%%%%%%%%%%%%%%%%

\section{Introduction}
\label{sec:intro}
CMC surfaces in space forms are given by their associated $\mathbb C^*$-family of flat $\SL(2,\C)$-connections
 $\nabla^\lambda$.
For simply connected surfaces, all flat connections are trivial, and the associated family of gauge equivalence classes $[\nabla^\lambda]$
does not contain any information. This is in stark contrast to compact surfaces of genus $g\geq1.$
For CMC tori, it was shown by Hitchin \cite{H} that the conformal type of the torus (as well as its energy) are determined by the asymptotic behavior of $[\nabla^\lambda]$ as $\lambda\to0.$ Similarly, it was shown in \cite{He3} that
the energy of Lawson symmetric CMC surfaces of genus $2$ can also be recovered from the asymptotic behavior of  $[\nabla^\lambda].$
In this paper, we prove that also the conformal type of a  compact  CMC surface of genus $g\geq2$ with simple umbilics is determined by the asymptotic behavior of the associated family of gauge equivalence classes $[\nabla^\lambda]$ locally in the Teichm\"uller space:
The trace of the monodromy of $\nabla^\lambda$ along curves $\gamma$ (satisfying a certain condition)
give the period of the square root of the Hopf differential along $\gamma,$ see formula \eqref{asymptotic_trace}. 
Then, we  reprove that the periods of the square root of a holomorphic quadratic differential with simple zeros (defined on the so-called Hitchin covering)
determine the Riemann surface structure locally in the Teichm\"uller space. We then apply this observation to
prove our main theorem.

\section{The family of flat connections and its asymptotic monodromy}\label{family}
Let $f\colon M \to S^3$ be a conformal immersion of constant mean curvature $H$  from a compact Riemann surface to the round $S^3.$ The immersion $f$ is determined by its associated
$\C^*$-family $\nabla^{\lambda}$ of flat $\SL(2, \C)$ connections on the trivial $\C^2$ bundle over $M$, see \cite{H, B}:

\begin{The}[\cite{H,B}]
\label{The1}
For a conformal CMC immersion $f\colon M\to S^3$  there exists
an associated family of flat $\SL(2,\C)$-connections
\begin{equation}\label{associated_family}
\lambda\in\C^*\mapsto \nabla^\lambda=\nabla+\lambda^{-1}\Phi-\lambda\Phi^*
\end{equation}
on a complex  bundle $V\to M$ of rank $2,$ where $\Phi$ is a nilpotent and nowhere vanishing
complex linear $1$-form and $\Phi^*$ is its
adjoint with respect to a unitary metric. The connections $\nabla^\lambda$
are unitary for $ S^1\subset\C^*$ and trivial for
$\lambda_1\neq\lambda_2\in  S^1$.  \\

Conversely, for every $\C^*$-family of flat connections satisfying the properties listed above,  the $SU(2)$-valued gauge transformation between the trivial  connections
$\nabla^{\lambda_1}$ and $\nabla^{\lambda_2}$ is a CMC immersion into $S^3=SU(2)$ with mean curvature $H=i\frac{\lambda_1+\lambda_2}{\lambda_1-\lambda_2}.$ \end{The}

Because of the following theorem  we cannot apply the methods developed by Hitchin \cite{H} in order to study
CMC surfaces of higher genus.

\begin{The}[\cite{He1}]
For compact CMC surfaces of genus $g\geq2$ the flat $SL(2, \C)$-connections $\nabla^{\lambda}$ are irreducible for generic spectral parameter $\lambda\in\C^*.$
\end{The}

\subsection{Monodromy}\label{mon}
Because $\Phi$ in Theorem \ref{The1} is nilpotent and nowhere vanishing, there exists a holomorphic line bundle $S^*=\ker \Phi\subset V$
which automatically satisfies $(S^*)^2=K^*,$ i.e., $S^*$ is the dual of a holomorphic spin bundle $S.$ With respect to
the unitary decomposition
$V=S^*\oplus S$
the connections $\nabla^\lambda$  are given by
\begin{equation}\label{gceqn}
\nabla^\lambda=\dvector{&\nabla^{spin^*} & 
-\frac{i}{2} Q^*\\ & -\frac{i}{2} Q & \nabla^{spin}}+\lambda^{-1}\dvector{0 & \tfrac{1}{2}\\0&0} -\lambda\dvector{0&0\\i\vol &0},\end{equation}
see \cite{He1,He3}.
In \eqref{gceqn} $Q$ is a holomorphic quadratic differential -- the Hopf differential of the CMC surface -- and $\nabla^S$ is the spin connection of the induced Levi-Civita connection. 
Let $\zeta^2=\lambda$ and consider the $\zeta$-dependent gauge transformations
\[g=\dvector{\tfrac{1}{\sqrt{\zeta}} &0\\0& \sqrt{\zeta}}\]
with respect to the decomposition $V=S^*\oplus S.$
Then, 
\begin{equation}\label{twisted}
\nabla^\zeta:=\nabla^\lambda.g=\tilde\nabla+\zeta^{-1}\dvector{0 & \tfrac{1}{2}\\-\frac{i}{2} Q&0} -\zeta\dvector{0&-\frac{i}{2} Q^*\\i\vol &0},
\end{equation}
is the twisted family of flat connections. For $\zeta^2=\lambda$ the traces of the monodromies of $\nabla^\zeta$ and $\nabla^\lambda$
along any closed curve are the same.

Instead of working with $\nabla^\lambda$ or $\nabla^\zeta$ we will work with another family of flat connections:
Consider the Hopf differential $Q$ of the CMC surface of genus $g\geq2$ and assume that it
has only simple zeros. Then there exists a double covering $\hat M\to M$ (the Hitchin covering \cite{H1}) branched over the umbilics, i.e., the zeros of $Q$.
$\hat M$ has genus $4g-3$ and there exists a holomorphic 1-form $\omega$ on $\hat M$ such that $\tfrac{i}{4}Q=-\omega^2$. Note that in \cite{H1} $Q$ is pulled back as a section and not as a quadratic differential. Nevertheless, it follows easily
that the pullback of $Q$ as a quadratic differential has a globally well-defined square root which is a holomorphic 1-form on the Hitchin covering.

 The holomorphic Higgs field $\Phi=\dvector{0 & \tfrac{1}{2}\\-\frac{i}{2} Q&0}$ (with respect to $\dbar^{\tilde\nabla}$) diagonalizes on $\hat M$ \cite{H1},
i.e., there is a decomposition $V=L\oplus L^*\to\hat M$ such that 
\[\Phi=\dvector{\omega&0\\0&-\omega},\]
where $\Phi$ is the pull-back as an endomorphism-valued 1-form.
In the following, we will work with the pull-back connections $\hat\nabla^\zeta$ 
of $\nabla^\zeta$ on $V=L\oplus L^*\to\hat M.$

In order to obtain
informations from the asymptotic behavior of the monodromy of $\hat\nabla^\zeta$ (respectively $\nabla^\lambda$)
we need to recall basic knowledge from the asymptotic analysis of linear ordinary differential equations. For a more extensive treatment of this topic we refer to \cite{Was,Fed}. A short treatment appropriate for our purposes is given in Appendix C in \cite{GMN}:
Let $I=(\varphi_0;\varphi_1)\subset\R$
and consider the sector \[S:=\{te^{i\varphi}\mid t\in\R^{>0},\varphi\in I\}.\]  Let $\gamma\colon [0;1]\to\hat M$ be a closed curve such that
\begin{equation}\label{direction_condition}
\Re(e^{i\varphi}\omega(\gamma'(t)))>0
\end{equation}
for all directions $e^{i\varphi}\in S^1$ with  $\varphi\in I$ and all
$t\in [0;1].$ Then the trace of the monodromy along $\gamma$ satisfies
\begin{equation}\label{asymptotic_trace}
\Tr(M(\hat\nabla^\zeta,\gamma))\exp{(-\tfrac{1}{\zeta}\int_\gamma\omega)}\xrightarrow[\tfrac{1}{\zeta}\in S,\,\zeta\to 0]{}  const\neq0\end{equation}
for some $const\in\C\setminus\{0\}.$ Hence, we obtain:
\begin{Pro}\label{as-be}
For closed curves $\gamma$ satisfying \eqref{direction_condition} the value of 
the period
$\int_\gamma\omega$
is determined by the asymptotic behavior for $\lambda\to0$  of the trace of the monodromy representation of $\nabla^\lambda.$
\end{Pro}

\section{The Riemann surface structure}\label{RSS}

Let $M_t,$ $t\in(-\delta,\delta),$ be a (smooth) family of compact Riemann surfaces of genus $g,$ and
let $Q_t$ be a family of holomorphic quadratic differentials  with respect to $M_t.$ 
Assume that the differentials $Q_t$ have only simple zeros, and consider
the family $\pi_t\colon\hat M_t\to M_t$ of Hitchin coverings as well as the corresponding (holomorphic) involutions 
$\sigma_t$ interchanging the sheets.
We can identify the smooth surfaces $\hat M=\hat M_t$ in such a way that all $\sigma_t$ are the same (smooth map $\sigma\colon \hat M\to\hat M$). In this situation we get:

\begin{Lem}\label{torelli}
Assume that
on the Hitchin covering $\hat M\to M$  the periods of $\omega_t=\sqrt Q_t\in H^0(\hat M_t,K_{\hat M_t})$  
are constant in $t.$ Then the Riemann surfaces $M_t$  are equivalent. Moreover, the Hopf differentials
$Q_t$ are also equivalent.
\end{Lem}
\begin{Rem}
It is a well-known fact (see \cite{KH} or $\S 1.1.$ in \cite{KZ}) that the space of Riemann surfaces equipped with 
holomorphic differentials $\omega$ having a prescribed number of zeros of prescribed order has local coordinates given 
by the so-called relative periods $\int_{\gamma_k}\omega,$ where the curves $\gamma_k,\, k=1,..,$ generate the first homology
group relative to the set of zeros of $\omega.$ 
We briefly reprove this fact in the situation of this paper along with
the necessary computation of  the ``half-periods'' $\int_\gamma\omega$ for curves connecting
zeros of $\omega.$
\end{Rem}
\begin{proof}
Consider the smooth surface $\hat M$ together with the family of closed complex valued 1-forms $\omega_t.$
They are holomorphic with respect to a $t$-dependent Riemann surface structure $J_t,$ and their zeros
are all of order $2.$ Hence, we can assume that the zeros coincide on the smooth surface $\hat M.$ Moreover,
we can also assume without loss of generality that in a local neighborhood $U\subset\hat M$ of its common zeros, the 1-forms coincide $(\omega_t)_{\mid U}=(\omega_0)_{\mid U}.$ 
As the periods of  $\omega_t$ are independent of $t$ we get
\[\omega_t=\omega_0+ df_t\]
for some complex-valued function $f_t\colon\hat M\to\C$ depending smoothly on $t.$
Moreover, the differential $df_t\in\Omega^1(\hat M,\C)$ vanishes  in the neighborhood $U$ of the zeros of $\omega.$ 
In order to apply the Moser trick to deform the forms $\omega_t$ to $\omega_0$, we first need to establish that the values of $\frac{\partial{f_t}}{\partial t}$  near the zeros of
$\omega_t$ are $0$.
%To prove this, it is enough to show that for all curves
%$\gamma$ connecting two zeros  we have
%\[\int_\gamma\omega_t=\int_\gamma\omega_0.\]
Recall that we have already  identified the involutions $\sigma_t=\sigma\colon \hat M\to\hat M.$
 Let $\gamma$ be a curve connecting
 two zeros of $\omega_t$ and consider the closed curve $\Gamma=\sigma(\gamma^{-1})\circ\gamma.$ 
 Because $\sigma^*\omega_t=-\omega_t$ for all $t$ we
 obtain
 \[2\int_\gamma\omega_t=\int_{\Gamma}\omega_t=\int_{\Gamma}\omega_0=2\int_\gamma\omega_0.\]
 Therefore, the values of $f_t$ at any zero of $\omega_t$ are the same, and without loss of generality we have that
 they are $0$ independently of $t.$
 Next we consider the closed $1$-form
 \[\Omega=\omega_0+df_t+\frac{\partial{f_t}}{\partial t} dt=\omega_t+\frac{\partial{f_t}}{\partial t} dt\]
 on $\hat M\times(-\delta,\delta).$ We want to show that there exists a (real) vector field $X$ on 
 $\hat M\times(-\delta,\delta)$ such that
 $\Omega(X)=0$ and $dt(X)=1.$ We already know that $f_t$ is $0$ in the neighborhood $U$ of the zeros of $\omega_t$
 for all $t\in(-\delta,\delta).$ Therefore, the existence and uniqueness of $X$ follows from the fact that
 for all $p\in \hat M\setminus \{\text{zeros of } \omega_t\}$ the linear map $(\omega_t)_p\colon T_p\hat M\to\C$ is an isomorphism. 
 Clearly, $X$ is complete.
 Then
 \begin{equation}\label{lieder}\mathcal L_X\Omega=d i_X\Omega+i_X d\Omega=0,\end{equation}
 and the flow at time $s$ \[\Phi_s\colon \hat M\times(-\delta,\delta)\to \hat M\times(-\delta,\delta)\]
 is of the form
 \[\Phi_s(p,t)=(\phi_s(p),t+s).\]
 By construction and \eqref{lieder}, $\phi_s\colon\hat M\to\hat M$ is a diffeomorphism which satisfies 
 $\phi_s^*\omega_s=\omega_0.$
 This also shows that $\phi_s$ is a holomorphic diffeomorphism with respect to the Riemann surface structures $M_0$ and $M_s.$
\end{proof}

\subsection{Deformations of CMC surfaces}
In the abelian case of tori, it is easy to show that the monodromy representation (based at some fixed point) already determines the conformal type of the immersed CMC torus $f\colon M\to S^3$:
The asymptotic behavior of the traces of the monodromies determines the periods of a (non-zero) holomorphic differential on the 
torus, see \cite{H} for details. Clearly, these periods determine the conformal type of the torus.
The main result of this note is the following generalization to compact CMC surfaces of genus $g\geq2:$

\begin{The}\label{Main}
Let $M$ be a compact oriented surface.
Let  $f_t\colon M\to S^3$  be a family of CMC surfaces
 such that 
the gauge equivalence classes of associated family of flat connections ${^t}{\nabla}^\lambda$
are equal for all $t$ and generic $\lambda\in\C^*.$ If the CMC surfaces have only simple umbilics (for generic $t$).
Then the induced Riemann surface structures $\Sigma_t$ are all equivalent and the Hopf differentials coincide. 
\end{The}
\begin{Rem}
The theorem remains true for CMC surfaces with periods, i.e., for CMC surfaces into a flat $S^3$ bundle over $M$
(in the sense of Hitchin's \cite{H} gauge-theoretic harmonic map equations). 
These are given by solutions $(g,Q)$ of the Gauss-Codazzi equations which are globally defined on the surface $M.$
Here $g$ is a Riemannian metric which induces the spin connection in\eqref{gceqn}, and $Q$ is the Hopf differential. The Gauss-Codazzi equations are then equivalent
to the flatness of all connections $\nabla^\lambda$ in \eqref{gceqn}. The theorem also remains true for CMC surfaces (with periods) in the space forms $\mathbb R^3$ and $\mathbb H^3.$

\end{Rem}
\begin{proof}
We  show that there exists closed curves $\gamma_i,$ $i=1,...,$ on $\hat M$
which satisfy \eqref{direction_condition} for non-empty sectors and which generate the first homology over $\mathbb Q.$
Then the proof follows from Lemma \ref{torelli} and Proposition \ref{as-be}.
We apply Lemma 4.4 in \cite{For}. This Lemma and its proof  tell us, that (the pull-back of) $Q$ on $\hat M$ is the limit of holomorphic quadratic differentials $Q_i$ with respect to suitable Riemann
surface structures $\hat M_i$ (note that we identify $\hat M_i$ and $\hat M$ as smooth surfaces) and which satisfy the following properties:
\begin{enumerate}
\item
all $Q_i$ have the same number of zeros of the same order as $Q;$ 
\item for all $i\in\N$ the horizontal distribution $\Re(\sqrt{Q_i}))=0$
contains closed curves which do not hit any zeros of $Q_i$ and which generate a Lagrangian subspace $\mathcal L$ in the symplectic space $H_1(\hat M,\R)$;
\item for any Lagrangian subspace $\Lambda\subset H_1(\hat M,\R)$ the differentials $Q_i$ can be chosen in such a way that $\mathcal L\cap \Lambda=\{0\};$
\item there is an open neighborhood $U$ of the zeros of $Q$ such that   $(Q_i)_{\mid U}=Q_{\mid U}$ for all $i\in\N.$ 
\end{enumerate}
We claim that for $Q_i$ close to $Q$ and a closed curve $\gamma\colon S^1\to\hat M\setminus\{\text{zeros of } Q_i\}$ satisfying
\[\Re(\sqrt{Q_i}(\gamma'))=0\]
condition \eqref{direction_condition} is satisfied for a non-empty sector $S$. This follows  from the fact that $\sqrt{Q_i}_{\mid U}=\sqrt{Q}_{\mid U}$ on $U$ and that point-wise  on the compact set $\hat M\setminus U$ the norms of $\sqrt{Q}$ and $\sqrt{Q_i}$ are bounded
from below uniformly. Note that because of (iii) we obtain enough curves $\gamma$ which generate the first homology $H_1(\hat M;\mathbb Q).$
\begin{Rem}
The statement of the theorem remains true for some classes of CMC surface with higher order umbilics. In particular, the conformal type of a $(k,l)$-symmetric CMC surface (as defined in \cite{HHSch}) is determined
uniquely by the asymptotic behavior \eqref{asymptotic_trace} of the monodromy representation.
\end{Rem}
\begin{Rem}
Combining Theorem \ref{Main} with Theorem 7 of \cite{He3} shows that deformations of compact CMC surfaces (with simple umbilics) which preserve the conjugacy classes of the monodromy representations of the associated families
of flat connections must be induced by a family of  dressing transformations. Dressing is a transformation of the family of flat connections induced by a $\lambda$-dependent gauge which becomes singular 
at certain spectral parameters $\lambda_0\notin S^1\cup\{0,\infty\}$ where parallel eigenlines with respect to $\nabla^{\lambda_0}$ exist. In the case of higher genus $g\geq2,$ families of dressing transformations
can only exists if there is a spectral parameter $\lambda_0$ for which $\nabla^{\lambda_0}$ is trivial (up to an diagonal $\Z_2$-gauge). We expect that this would occur only in very exceptional situations.
\end{Rem}

%The general case follows by approximating $f$ by conformal harmonic sections of a flat $S^3$ bundle over $M$
%(in the sense of Hitchin's \cite{H} gauge-theoretic harmonic map equations) which have 
%Hopf differentials with simple zeros.
\end{proof}

%\subsection{Dressing Transformations}

%%%%%%%%%%%%%%%%%%%%%%%%%%%%%%%%%%%%%%%%%%%%%%%%%%%%%%%%%%%%%%%%%%%%%%%%%%%%%%
%           Bibliography                                                     %
%%%%%%%%%%%%%%%%%%%%%%%%%%%%%%%%%%%%%%%%%%%%%%%%%%%%%%%%%%%%%%%%%%%%%%%%%%%%%%

\end{document}